\RequirePackage{fix-cm}
 \documentclass{article}
\usepackage{graphicx}

\usepackage{amsmath,amssymb}

\newcommand{\E}{\mathbb E}

\newcommand{\dd}{\text{\rm d}}
\newcommand{\x}{\mathbf x}
\newtheorem{theorem}{Theorem}
\newtheorem{lemma}{Lemma}
\newtheorem{corollary}{Corollary}
\newtheorem{remark}{Remark}
\begin{document}

\title{Central limit theorem for Gibbsian U-statistics of facet processes 
}


\author{Jakub Ve\v{c}e\v{r}a \\
\small{Charles University in Prague}, \\
\small{Department of Probability and Mathematical Statistics},\\
\small{vecera@karlin.mff.cuni.cz}
 }
\maketitle

\begin{abstract}
Special case of a Gibbsian facet process on a fixed window with a discrete orientation distribution and with increasing intensity of the underlying Poisson process is studied. All asymptotic joint moments for interaction U-statistics are calculated and using the method of moments the central limit theorem is derived.

\end{abstract}

\section{Introduction}

In the present paper we use methods developed in \cite{RefV} to calculate all moments of Gibbsian $U$-statistics of facets in a bounded window in arbitrary Euclidean dimension. These moments are used to derive the central limit theorem for such statistics. Central limit theorems for $U$-statistics of Poisson processes were derived based on Malliavin calculus and the Stein method in \cite{RefR}.

Our calculations are based on the achievements in \cite{RefB}, where functionals of spatial point processes given by a density with respect to the Poisson process were investigated using the Fock space representation from \cite{RefL}. This formula is applied to the product of a functional and the density and using a special class of functionals called $U$-statistics closed formulas for mixed moments of functionals are obtained. In processes with densities the key characteristic is the correlation function \cite{RefG} of arbitrary order which is dual to kernel function of the density as a function of the Poisson process. 
 
As in \cite{RefV} we call facets some compact subsets of hyperplanes with a given shape, size and orientation. Natural geometrical characteristics of the union of the facets, based on Hausdorff measure of the intersections of pairs, triplets, etc., of facets form $U$-statistics. Building a parametric density from exponential family, the limitations for the space of parameters have to be given, so called submodels are investigated. In application of the moment formulas we are interested in the limit behaviour when the intensity of the reference Poisson process tends to infinity.

 We restrict ourselves to the facet model with finitely many orientations corresponding to canonical vectors. 
 In \cite{RefV} basic asymptotic properties of the studied $U$-statistics are derived. When the order of the submodel is not greater than the order of the observed $U$-statistic then asymptoticaly the mean value of the $U$-statistic vanishes. This leads to a degeneracy in the sense that some orientations are missing. On the other hand  when the order of the submodel is greater than the order of the observed $U$-statistic then the limit of correlation function is finite and nonzero and under selected standardization $U$-statistic tends almost surely to its non-zero expectation. By changing the standardization, however, we achieve a finite non-zero asymptotic variance. In the present paper we simplify the calculation of moments so that we are able to calculate any asymptotic moment.
 
\section{Central limit theorem}
Let $Y=[0,b]^{d} \times \{2b \} \times \{ e_{1},\ldots,e_{d} \}$ be a space of facets (facets are $d-1$ dimensional cubes) with three parts: set of facet centres, possible sizes of facet and possible orientations of hyperplane containing facet, i.e. we consider only facet with fixed size and orientation described only by elementary vectors. Facet in such case can be described as
\begin{align}
\nonumber  ((z_{1},\ldots,z_{d}),2b,e_{l}) \rightarrow \{(x_{1},\ldots,x_{d}) \in \mathbb{R}^{d}, x_{l}=z_{l};|x_{i}-z_{i}|\leq b,i \in [d] \setminus \{l \} \},\\
\nonumber (z_{1},\ldots,z_{d}) \in [0,b]^{d}, l\in [d], [d]= \{1,\ldots,d \}
\end{align} 
Moreover let  $({\mathbf N},{\mathcal N})$ be a measurable space of integer-valued finite measures on $Y$, where ${\mathcal N}$ is the smallest $\sigma$-algebra which makes the mappings $\x\mapsto \x(A)$ measurable for all Borel sets $A\subset Y$ and all $\x\in {\mathbf N}$.  We denote $\eta_{a}$ finite Poisson process of facets with intensity function $a \lambda$ on $Y$, $a \geq 1$ in form
\begin{align}
\nonumber a\lambda (\dd x)=a\lambda (\dd (z,r,\phi)) = a\chi (z) \dd z \delta_{2b}(r) \frac{1}{d} \sum_{i=1}^{d} \delta_{e_{i}} (\dd \phi),
\end{align}
where we have fixed the facet size, uniform distribution of the facet orientation and $\chi:[0,b]^{d} \rightarrow \mathbb{R}_{+}$. We also define interaction $U$-statistics (using Hausdorff measure $\mathbb{H}^{d-j}$ of order $d-j$)
\begin{align}
\nonumber G_{j}(\mathbf{x}) = \frac{1}{j!} \sum_{(x_{1},\ldots,x_{j}) \in \mathbf{x}^{j}_{\neq}} \mathbb{H}^{d-j} (\cap_{i=1}^{j} x_{i} )
\end{align}
and the process $\mu_{a}$ with density 
\begin{align}
\nonumber p(\mathbf{x})=c_{a} \exp \left( \sum_{i=1}^{d} \nu_{i} G_{i}(\mathbf{x}) \right)
\end{align}
 with respect to $\eta_{a}$, where $a \geq 1$, $\nu_{i}$ is a real parameter and $ c_{a}=\frac{1}{\E \exp \left( \sum_{i=1}^{d} \nu_{i} G_{i}(\mathbf{x}) \right)}$. Fullfilling of condition $\nu_{i} \leq 0 ,i=2,\ldots,d$ assures that $p(\eta_{a}) \in L^{1}(P_{\eta_{a}}) \cap L^{2}(P_{\eta_{a}}) $, we will discuss necessity of these conditions later. We also use the notion of a submodel $\mu_{a}^{(l)}$, where $\nu_{j}=0, j\neq l$ and $\nu_{l} < 0$. We will explore properties of such submodels with the order higher than $1$, because in case $\mu_{a}^{(1)}$ we deal with Poisson process \cite{RefV}. \\

 We can use a short expression for moment formulas using diagrams and partitions, see \cite{RefP}, \cite{RefL}. Let $\tilde{\prod}_k$ be the set of all  partitions $\{J_i\}$ of $[k],$ where $J_i$ are disjoint blocks and $\cup J_i=[k].$ For
$k=k_1+\dots +k_m$ and blocks $$J_i=\{ j: k_1+\dots +k_{i-1}< j\leq k_1+\dots +k_i\},\; i=1,\dots ,m,$$ consider the partition $\pi =\{J_i,\;1\leq i\leq m\}$ and let
$\prod_{k_1,\dots ,k_m}\subset \tilde{\prod}_k$ be the set of all partitions $\sigma\in\tilde{\prod}_k$ such that $|J\cap J'|\leq 1$ for all $J\in\pi $ and all $J'\in\sigma .$ Here $|J|$ is the cardinality of a block $J\in\sigma .$ We will be referring to blocks of $\pi$ as to rows and we denote $S(\sigma)$ the number of pure singleton rows of partition $\sigma$, i.e. number of $J \in \pi$, which for all $J' \in \sigma$, $|J \cap J'|=1 \Rightarrow |J'|=1$. \\

For a partition $\sigma\in\prod_{k_1\dots k_m}$ and measurable functions $f_{i}:B \rightarrow \mathbb{R},j=1,\ldots,m,$ we define the function $(\otimes_{j=1}^mf_j)_\sigma:B^{|\sigma |}\rightarrow {\mathbb R}$ by replacing all variables of the tensor product $\otimes_{j=1}^mf_j$ that belong to the same block of $\sigma $ by a new common variable, $|\sigma |$ is the number of blocks in $\sigma .$ We denote $\Pi^{(m_{1},\ldots,m_{s}) }_{1,\ldots,s}=\Pi_{1,\ldots,1,\ldots,s,\ldots,s}$, where $i$ repeats $m_{i}$ times for $i=1,\ldots,s$. It holds 
$\binom{p}{q} =0,q>p$. Now we state the main theorem of the paper.

\begin{theorem}\label{CLVT}
Denote $\tilde{G}_{j}(\mu^{(c)}_{a}) = \frac{ G_{j}(\mu^{(c)}_{a})- \E G_{j}(\mu^{(c)}_{a})}{a^{j-\frac{1}{2}}}$, $ 1 \leq j  \leq d$, $2 \leq c \leq d$, then
\begin{align}\label{CLV}
(\tilde{G}_{1}(\mu^{(c)}_{a}),\ldots,\tilde{G}_{d}(\mu^{(c)}_{a})) \xrightarrow{\mathcal{D}} \mathbf{Z},~c=2,\ldots,d
\end{align}
as $a$ tends to infinity, where $\mathbf{Z} \sim N(0,\Sigma)$, $\Sigma =\{ \theta_{ij} \}_{i,j=1}^{d}$,
\begin{align}
\nonumber \theta_{kl} &=\frac{(c-1)}{d^{k+l-1}} \binom{c-2}{k-1} \binom{c-2}{l-1}I_{kl},\\ 
\nonumber I_{kl} &= \int_{([0,b]^{d})^{k+l-1}} {\mathbb H}^{d-k} (\cap_{i=1}^{k} (s_{i},2b,e_{i})) {\mathbb H}^{d-l} (\cap_{i=2}^{l} (s_{i+k-1},2b,e_{i}) \cap (s_{1},2b,e_{1}) ) \times \\
\nonumber & \times \chi(s_{1}) \dd s_{1},\ldots,\chi(s_{k+l-1}) \dd s_{k+l-1},
\end{align}
moreover 
\begin{align}
\nonumber &G_{j}(\mu^{(c)}_{a}) \xrightarrow{L^{2}} 0 ~, &  c \in \{ 2,\ldots,d \}, j\geq c  , \\
\nonumber &\frac{G_{j}(\mu^{(c)}_{a})}{a^{j}} \xrightarrow{L^{2}} \frac{I_{j}}{d^{j}} \binom{c-1}{j} ~, & c \in \{ 2,\ldots,d \}, j < c ,
\end{align}
where $I_{j}= \int_{([0,b]^{d})^{j}} {\mathbb H}^{d-j}(\cap_{i=1}^{j} (s_{i},2b,e_{i}) )  \chi(s_{1}) \dd s_{1},\ldots,\chi(s_{j}) \dd s_{j}$.

\end{theorem}
\begin{remark}
Random variables $\tilde{G}_{c}(\mu^{(c)}_{a}),\tilde{G}_{c+1}(\mu^{(c)}_{a}),\ldots,\tilde{G}_{d}(\mu^{(c)}_{a})$ are asymptotically degenerate, i.e. their expectations tend to zero, thus covariance of these variables $\theta_{kl}=0,k \geq c,l \in [d]$.
\end{remark}

\begin{remark}
For random vector $(\tilde{G}_{1}(\eta_{a}),\ldots,\tilde{G}_{d}(\eta_{a}))$ we have similar results \cite{RefLST} with $\theta_{kl}=\frac{d}{d^{k+l-1}} \binom{d-1}{k-1}\binom{d-1}{l-1}I_{kl}$.
\end{remark}

\begin{corollary}
It holds
\begin{align}
\frac{ G_{j}(\mu^{(c)}_{a})- \E G_{j}(\mu^{(c)}_{a})}{a^{j-\frac{1}{2}}}  \xrightarrow{\mathcal{D}} Z, c=2,\ldots,d,j<c
\end{align}
as $a$ tends to infinity, where  $Z \sim N(0,\theta_{jj})$.

\end{corollary}

\begin{lemma}\label{KorelFce}
It holds
\begin{align}
\label{KorelLimita}\rho_{p}(x_{1},\ldots,x_{p},\mu^{(c)}_{a})= \frac{\E \exp(\nu_{c} G_{c}(\eta_{a}\cup \{x_{1},\ldots,x_{p}\} ))}{\E \exp(\nu_{c} G_{c}(\eta_{a} ))} \rightarrow  
\frac{\frac{(d-k)!}{(c-1-k)!}}{\frac{d!}{(c-1)!}},
\end{align}
as $a$ tends to infinity, where $x_{i} \in Y$ and $k$ is number of facet orientations among $\{x_{1},\ldots,x_{p} \}$ and $c \geq 2$. Moreover 
\begin{align}
\nonumber \left\vert \rho_{p}(x_{1},\ldots,x_{p},\mu_{a}^{(c)}) - \frac{\frac{(d-k)!}{(c-1-k)!}}{\frac{d!}{(c-1)!}} \right\vert < Re^{-Sa},R,S>0,
\end{align}
where $R$ and $S$ do not depend on $x_{1},\ldots,x_{p}$.
\end{lemma}
\begin{remark}
The function $\rho_{p}$ is called the correlation function  of order p \cite{RefG}.
\end{remark}
\noindent{\bf Proof:}
First consider submodel $\mu_{a}^{(c)}$ and facets $x_{1},\ldots,x_{p}$ with $ p \leq c$ distinct orientations and without loss of generality consider orientations $e_{1},\ldots,e_{p}$, because the distribution of orientations is uniform and the $\rho_{p}(x_{1},\ldots,x_{p},\mu^{(c)}_{a})$ does not change under rotation uniformly applied to all facets $x_{1},\ldots,x_{p}$.
It holds $\rho_{p}(x_{1},\ldots,x_{p},\mu^{(c)}_{a})=$
\begin{align}
\nonumber  =\frac{\sum_{n=0}^{\infty}\frac{a^{n}}{n!} \int_{Y^{n}} \exp\left(    \nu_{c}G_{c} \{ u_{1},\ldots,u_{n},x_{1},\ldots,x_{p} \} \right) \lambda^{n} (\dd (u_{1},\ldots,u_{n}))}{\sum_{n=0}^{\infty}\frac{a^{n}}{n!} \int_{Y^{n}} \exp\left(    \nu_{c}G_{c} \{ u_{1},\ldots,u_{n} \} \right) \lambda^{n} (\dd (u_{1},\ldots,u_{n}))}.
\end{align}
We can obtain bounds for this expression by using the bounds for the volumes of intersection of facets $b^{d-c} \leq \mathbb{H}^{d-c}(\cap_{i=1}^{c}y_{i}) \leq (2b)^{d-c} $ as follows
\begin{gather}
\nonumber \frac{\sum_{n=0}^{\infty}\frac{(\frac{aT}{d})^{n}}{n!} \sum_{n_{1}+ \ldots + n_{d}=n} \binom{n}{n_{1},\ldots,n_{d}} \exp \left( \nu_{c} (2b)^{d-c} R^{c,p}(p,d,\mathbf{n}^{(d)})\right)}{\sum_{n=0}^{\infty}\frac{(\frac{aT}{d})^{n}}{n!} \sum_{n_{1}+ \ldots + n_{d}=n} \binom{n}{n_{1},\ldots,n_{d}} \exp \left( \nu_{c} b^{d-c} R^{c,0}(0,d,\mathbf{n}^{(d)})\right)} \leq \\
\label{RHObounds}  \rho_{p}(x_{1},\ldots,x_{p},\mu^{(c)}_{a}) \leq \\
\nonumber  \frac{\sum_{n=0}^{\infty}\frac{(\frac{aT}{d})^{n}}{n!} \sum_{n_{1}+ \ldots + n_{d}=n} \binom{n}{n_{1},\ldots,n_{d}} \exp \left( \nu_{c} b^{d-c} R^{c,p}(p,d,\mathbf{n}^{(d)})\right)}{\sum_{n=0}^{\infty}\frac{(\frac{aT}{d})^{n}}{n!} \sum_{n_{1}+ \ldots + n_{d}=n} \binom{n}{n_{1},\ldots,n_{d}} \exp \left( \nu_{c} (2b)^{d-c} R^{c,0}(0,d,\mathbf{n}^{(d)})\right)},
\end{gather}
where $T= \int_{[0,b]^{d}} \chi(z) \dd z$, $n_{i}$ are the numbers of facets among $u_{1},\ldots,u_{n}$ with orientations $e_{i}$, respectively, $i=1,\ldots,d$ and $\mathbf{n}^{(d)}=(n_{1},\ldots,n_{d})$. Furthermore denote 
\begin{align}
\nonumber R^{c,p}(q,d,\mathbf{n}^{(d)}) = \sum_{\substack{F \subset [d] \\
c-p \leq |F| \leq c \\
|F \cup [q]|+p-q \geq c
}} \prod_{j \in F } n_{j}.
\end{align}
Specially $R^{c,0}(0,d,\mathbf{n}^{(d)})$ is the total number of intersections of all $c$-tuples of the facets among $u_{1},\ldots,u_{n}$ and $R^{c,p}(p,d,\mathbf{n}^{(d)})$ is the total number of intersections of all $c$-tuples of the facets among facets $u_{1},\ldots,u_{n},x_{1},\ldots,x_{p}$.
Then we substitute $\frac{aT}{d}$ for $a$, extend the both fractions by $e^{-a(c-1)}$ and then we get in the case of the lower bound  of (\ref{RHObounds})
\begin{align}
\nonumber \frac{\sum_{n_{1}=0}^{\infty} \ldots \sum_{n_{d}=0}^{\infty} \frac{a^{n_{1} + \ldots \ n_{d}}}{n_{1}! \ldots n_{d}!} \exp \left( \nu_{c} (2b)^{d-c} R^{c,p}(p,d,\mathbf{n}^{(d)}) - a(c-1)\right)}{\sum_{n_{1}=0}^{\infty} \ldots \sum_{n_{d}=0}^{\infty} \frac{a^{n_{1} + \ldots \ n_{d}}}{n_{1}! \ldots n_{d}!} \exp \left( \nu_{c} b^{d-c} R^{c,0}(0,d,\mathbf{n}^{(d)}) - a(c-1)\right)}
\end{align}
and in the case of the upper bound of (\ref{RHObounds}) 
\begin{align}
\nonumber \frac{\sum_{n_{1}=0}^{\infty} \ldots \sum_{n_{d}=0}^{\infty} \frac{a^{n_{1} + \ldots \ n_{d}}}{n_{1}! \ldots n_{d}!} \exp \left( \nu_{c} b^{d-c} R^{c,p}(p,d,\mathbf{n}^{(d)}) - a(c-1)\right)}{\sum_{n_{1}=0}^{\infty} \ldots \sum_{n_{d}=0}^{\infty} \frac{a^{n_{1} + \ldots \ n_{d}}}{n_{1}! \ldots n_{d}!} \exp \left( \nu_{c} (2b)^{d-c} R^{c,0}(0,d,\mathbf{n}^{(d)}) - a(c-1)\right)}.
\end{align}
Using Lemma \ref{LemmaLim} we get the limit of the lower and upper bound in the same form  $ \frac{\frac{(d-p)!}{(c-1-p)!}}{\frac{d!}{(c-1)!}} $. For $d \geq p>c$  we can get an upper bound in (\ref{RHObounds}) just by using $p=c$, which tends to zero. \\
Now consider more than one facet with the same orientation among $x_{1},\ldots x_{p}$ and with $k<c$ distinct orientations, which are without loss of generality set to $e_{1},\ldots,e_{k}$ and maximum number of facets with the same orientation is $P$, then we can bound the correlation function

\begin{gather}
\nonumber \frac{\sum_{n=0}^{\infty}\frac{(\frac{aT}{d})^{n}}{n!} \sum_{n_{1}+ \ldots + n_{d}=n} \binom{n}{n_{1},\ldots,n_{d}} \exp \left( \nu_{c} P^{d} (2b)^{d-c} R^{c,p}(p,d,\mathbf{n}^{(d)})\right)}{\sum_{n=0}^{\infty}\frac{(\frac{aT}{d})^{n}}{n!} \sum_{n_{1}+ \ldots + n_{d}=n} \binom{n}{n_{1},\ldots,n_{d}} \exp \left( \nu_{c} b^{d-c} R^{c,0}(0,d,\mathbf{n}^{(d)})\right)} \leq \\
\nonumber  \rho_{p}(x_{1},\ldots,x_{p},\mu^{(c)}_{a}) \leq \\
\nonumber  \frac{\sum_{n=0}^{\infty}\frac{(\frac{aT}{d})^{n}}{n!} \sum_{n_{1}+ \ldots + n_{d}=n} \binom{n}{n_{1},\ldots,n_{d}} \exp \left( \nu_{c} b^{d-c} R^{c,p}(p,d,\mathbf{n}^{(d)})\right)}{\sum_{n=0}^{\infty}\frac{(\frac{aT}{d})^{n}}{n!} \sum_{n_{1}+ \ldots + n_{d}=n} \binom{n}{n_{1},\ldots,n_{d}} \exp \left( \nu_{c} P^{d} (2b)^{d-c} R^{c,0}(0,d,\mathbf{n}^{(d)})\right)}.
\end{gather}

These bounds lead to expressions in the same form as in the case with different orientations and thus we proceed in the same way and get the value of the limit $\frac{\frac{(d-k)!}{(c-1-k)!}}{\frac{d!}{(c-1)!}}$. For $d \geq k \geq c$ we need only lower bound for the number of intersections in form $R^{c,k}(k,d,\mathbf{n}^{(d)})$, which forms upper bound for the correlation function, which tends to zero.\\
Bounds for the numerator and denominator of the correlation function converge to their limits with at least exponential rate and we can also see that the upper bounds can be selected to depend only on the $\nu$, $s$ and $c$ , therefore they do not depend on currently selected facets $x_{1},\ldots,x_{p}$ in the argument of correlation function. The rate of convergence can be extended to the whole fraction in \ref{KorelLimita} following way
\begin{align}
\nonumber |A(a)-A|< R_{1}e^{-S_{1}a},|B(a)-B|< R_{2}e^{-S_{2}a} \\
\nonumber \frac{A(a)}{B(a)} - \frac{A}{B} \leq \frac{\frac{A(a)}{1-R_{2}e^{-S_{2}a}}-A}{B} \leq \frac{R_{1}e^{-S_{1}a}+ A R_{2}e^{-S_{2}a}}{(1-R_{2}e^{-S_{2}a})B} \leq R e^{-Sa} \\
\nonumber \frac{A(a)}{B(a)} - \frac{A}{B} \geq  \frac{\frac{A(a)}{1+R_{2}e^{-S_{2}a}}-A}{B} \geq \frac{R_{1}e^{-S_{1}a}- A R_{2}e^{-S_{2}a}}{(1+R_{2}e^{-S_{2}a})B} \geq - R e^{-Sa},
\end{align}
where $A(a)$ is the value of numerator and $B(a)$ is the value of denominator on the left side in  (\ref{KorelLimita}),respectively $A$, $B$ is the limit of the numerator, denominator on the right side in (\ref{KorelLimita}).

\hfill $\Box $

\begin{lemma}\label{LemmaLim}
\begin{align}
\nonumber  \left\vert \sum_{n_{1}=0}^{\infty} \ldots \sum_{n_{d}=0}^{\infty} \frac{a^{n_{1}+\ldots + n_{d}}}{n_{1}! \ldots n_{d}!} \exp \left( \nu R^{c,p}(p,d,\mathbf{n}^{(d)})  - a(c-1) \right) - \frac{(d-p)!}{(c-1-p)!} \right\vert <R e^{-Sa}, 
\end{align}
where $\nu <0$, $p \leq c \leq d$, $c \geq 2$ and $R,S>0$.
\end{lemma}

\noindent{\bf Proof:}
We denote 
\begin{align}
\nonumber I(a,c,t,s)= \sum_{n_{1}=0}^{\infty} \ldots \sum_{n_{s}=0}^{\infty} \frac{a^{n_{1}+\ldots + n_{s}}}{n_{1}! \ldots n_{s}!} \exp \left( \nu R^{c,p}(t,s,\mathbf{n}^{(s)})  - a(c-1) \right).
\end{align}
Firstly we calculate the values of the limit by calculating the sum over \\$\left(n_{1}>0 \wedge \ldots \wedge n_{d}>0 \right)$ to show that this value tends to zero as $a$ tends to infinity. We show this only for $p=0$ because for $p>0$, we get upper bound using $p=0$ and the sum is non-negative.
In the following we use Chernoff bound for tail probabilities of Poisson distribution
\begin{align}
\nonumber \sum_{l=0}^{m} \frac{s^{l}}{l!} \leq \frac{(es)^{m}}{m^{m}}, m < s.
\end{align}
\begin{enumerate}

\item First we consider that all the summing variables are between $0$ and $a^{2/3}$: $\left( (a^{2/3} >n_{1}>0) \wedge \ldots \wedge (a^{2/3}>n_{d}>0) \right)$
\begin{align}
\nonumber  \sum_{n_{1}>0}^{a^{2/3}}\ldots\sum_{n_{d}>0}^{a^{2/3}} \frac{a^{n_{1}+\ldots+n_{d}}}{n_{1}! \ldots n_{d}!} \exp \left(\nu  R^{c,0}(0,d,\mathbf{n}^{(d)}) -a(c-1)\right) \leq \\
\nonumber  \sum_{n_{1}>0}^{a^{2/3}}\ldots\sum_{n_{d}>0}^{a^{2/3}} \frac{a^{n_{1}+\ldots+n_{d}}}{n_{1}! \ldots n_{d}!} \exp \left(-a(c-1)\right)  \leq \left( \frac{(ea)^{d a^{2/3}}}{(a^{2/3})^{d a^{2/3}}} \right)e^{-a(c-1)}\rightarrow 0,
\end{align}
where we used $d$-times the Chernoff bound.
\item Now consider that one of the summing variables is greater than $a^{2/3}$, without loss of generality we select $n_{d}$\\  $\left( (a^{2/3} >n_{1}>0) \wedge \ldots \wedge (a^{2/3}>n_{d-1}>0) \wedge (n_{d}\geq a^{2/3}) \right)$
\begin{align}
\nonumber   \sum_{n_{1}>0}^{a^{2/3}} \ldots\sum_{n_{d-1}>0}^{a^{2/3}}\sum_{n_{d}\geq a^{2/3}}^{\infty} \frac{a^{n_{1}+\ldots+n_{d}}}{n_{1}! \ldots n_{d}!} \exp \left( \nu  R^{c,0}(0,d,\mathbf{n}^{(d)}) - a(c-1) \right)  \leq \\
\nonumber \sum_{n_{1}>0}^{a^{2/3}} \ldots\sum_{n_{d-1}>0}^{a^{2/3}}\sum_{n_{d}=0}^{\infty} \frac{a^{n_{1}+\ldots+n_{d}}}{n_{1}! \ldots n_{d}!} \exp \left( \nu  n_{d} - a(c-1) \right)=\\
\nonumber \sum_{n_{1}>0}^{a^{2/3}} \ldots\sum_{n_{d-1}>0}^{a^{2/3}}  \frac{a^{n_{1}+\ldots+n_{d-1}}}{n_{1}! \ldots n_{d-1}!} \exp \left(ae^{\nu} -a(c-1)  \right) \leq \\
\nonumber \left( \frac{(ea)^{(d-1) a^{2/3}}}{(a^{2/3})^{(d-1) a^{2/3}}} \right)\exp \left(ae^{\nu} -a(c-1)  \right) \rightarrow 0
\end{align}
because  $e^{\nu}-(c-1)<0$.
\item When at least two of the summing variables are greater than $a^{2/3}$, without loss of generality select $n_{d-1}$ and $n_{d}$, then we have
\begin{align}
\nonumber  & \sum_{n_{1}>0}^{a^{2/3}} \ldots\sum_{n_{d-2}>0}^{a^{2/3}}\sum_{n_{d-1}\geq a^{2/3}}^{\infty}\sum_{n_{d} \geq a^{2/3}}^{\infty} \frac{a^{n_{1}+\ldots+n_{d}}}{n_{1}! \ldots n_{d}!} \exp \left(\nu  R^{c,0}(0,d,\mathbf{n}^{(d)}) -a(c-1)  \right)  \\
\nonumber & \leq \sum_{n_{1}>0}^{a^{2/3}} \ldots\sum_{n_{d-2}>0}^{a^{2/3}}\sum_{n_{d-1}\geq a^{2/3}}^{\infty}\sum_{n_{d} \geq  a^{2/3}}^{\infty} \frac{a^{n_{1}+\ldots+n_{d}}}{n_{1}! \ldots n_{d}!} \exp \left(\nu a^{4/3} -a(c-1) \right)  \\
\nonumber & \leq \exp \left( \nu a^{4/3} +a(d+1-c) \right) \rightarrow 0.
 \end{align}
\item The same applies to the case, where more than two variables are greater than $a^{2/3}$, because we are able to find terms with higher power of $a$ in the exponential.
\end{enumerate}
 Therefore we need only to examine the remaining terms, where at least one of the variables is equal to zero, thus we replace $I(a,c,p,d)$ by $d$ sums, where one variable is set to zero
 \begin{align}\label{IRekVzorec}
  I(a,c,p,d)\approx p I(a,c,p-1,d-1)+ (d-p)I(a,c,p,d-1), 
 \end{align}
where $\approx$ is the equality after omiting the summands, which tend to zero on the left side, $I(a,c,p-1,d-1)$ is the sum after setting to zero one of the variables $n_{1},\ldots,n_{p}$, $I(a,c,p,d-1)$ is the sum after setting to zero one of the variables $n_{p+1},\ldots,n_{d}$ and the multiplying numbers are the counts of possible selections of these variables. It can be shown that
\begin{align}
\nonumber &\lim_{a\rightarrow \infty}I(a,c,t,c-1) &= 0, & ~ t<p, \\
\nonumber & &=1, &  ~ t=p.
\end{align} 
Because the series on the right side of (\ref{IRekVzorec}) are in the same form as the original one and we can again sum only over the indices, where at least one is equal to zero, thus we repeat the $(d-c+1)$-times step in (\ref{IRekVzorec}) and we get 
\begin{align}
 I(a,c,p,d)\approx \sum_{j=0}^{d-c+1} c_{j} I(a,c,p-j,c-1),
\end{align}
where $c_{j} \in \mathbb{N}$. All summands tend to zero with one exception of \\$c_{0}I(a,c,p,c-1)$ with 
\begin{align}
\nonumber c_{0} &= \frac{(d-p)!}{(c-1-p)!}, c  > p, \\
\nonumber  &= 0, c  = p,
\end{align}
which is the number of possible selections of variables set to zero from $n_{p+1},\ldots,n_{d}$ in $d-c+1$ steps. The overall speed of convergence is implied by the convergence speed of every part of the sum, which converges to its limit at least at exponential rate.
\hfill $\Box $

\begin{lemma}\label{SimpleLemma}
It holds
\begin{align}
\nonumber \bigg| \int_{Y^{|\sigma|}}  \left( \otimes_{j=0}^{d-1}  \left(  \left( \bar{\mathbb{H}}^{d-j} \right)^{\otimes m_{j}}  \right) \right)_{\sigma} (u_{1},\ldots,u_{|\sigma|}) \rho_{|\sigma|} (u_{1},\ldots,u_{|\sigma|},\mu^{(c)}_{a}) \lambda^{|\sigma|} (\dd (u_{1},\ldots,u_{|\sigma|})) - \\
 \nonumber   \int_{Y_{c-1}^{|\sigma|}}  \left( \otimes_{j=0}^{d-1}  \left( \left( \bar{\mathbb{H}}^{d-j} \right)^{\otimes m_{j}}  \right) \right)_{\sigma}  (u_{1},\ldots,u_{|\sigma|})  \lambda^{|\sigma|} (\dd (u_{1},\ldots,u_{|\sigma|})) \bigg|< R e^{-Sa} ,
\end{align}
where $\sigma \in \Pi_{1,\ldots,s}^{(m_{1},\ldots,m_{s})}$, $S,R>0$, $Y_{c-1}= [0,b]^{d} \times \{ 2b \} \times \{ e_{1},\ldots, e_{c-1} \}$ is space of facets with  $d-c+1$ orientations missing (we can select any orientations), $m_{j} \in \mathbb{N}^{0}$ and $\bar{\mathbb{H}}^{d-j} (u_{1},\ldots,u_{j})= \mathbb{H}^{d-j} (\cap_{i=1}^{j} u_{i})$.
\end{lemma}
\begin{remark}
In space $Y_{c-1}$ we still use intensity $\lambda$, even when the orientations are restricted to $\{e_{1},\ldots,e_{c-1} \}$ and the measure is non-zero on orientations $\{e_{1},\ldots,e_{d} \}$.
\end{remark}

\noindent{\bf Proof:}

The limit of correlation function depends only on the number $k$ of the distinct orientations among the facets $(u_{1},\ldots,u_{|\sigma|})$, then correlation function tends to $\frac{ \frac{(d-k)!}{ (c-1-k)!}}{\frac{d!}{(c-1)!}}$ and thus we can write
\begin{gather}
\nonumber \int_{Y^{|\sigma|}} \left( \otimes_{j=0}^{d-1}  \left(\left( \bar{\mathbb{H}}^{d-j} \right)^{\otimes m_{j}}  \right) \right)_{\sigma} (u_{1},\ldots,u_{|\sigma|}) \rho_{|\sigma|} (u_{1},\ldots,u_{|\sigma|},\mu^{(c)}_{a}) \lambda^{|\sigma|} (\dd (u_{1},\ldots,u_{|\sigma|})) = \\
\nonumber \sum_{k=1}^{d} \binom{d}{k} \int_{(Y^{|\sigma|})_{[k]}} \left( \otimes_{j=0}^{d-1}  \left(\left( \bar{\mathbb{H}}^{d-j} \right)^{\otimes m_{j}}  \right) \right)_{\sigma} (u_{1},\ldots,u_{|\sigma|}) \times \\
\nonumber \times \rho_{|\sigma|} (u_{1},\ldots,u_{|\sigma|},\mu^{(c)}_{a}) \lambda^{|\sigma|} (\dd (u_{1},\ldots,u_{|\sigma|})) \leq  \\
\nonumber  \sum_{k=1}^{d} \binom{d}{k} \int_{(Y^{|\sigma|})_{[k]}} \left( \otimes_{j=0}^{d-1}  \left(\left( \bar{\mathbb{H}}^{d-j} \right)^{\otimes m_{j}}  \right) \right)_{\sigma} (u_{1},\ldots,u_{|\sigma|}) \frac{ \frac{(d-k)!}{ (c-1-k)!}}{\frac{d!}{(c-1)!}} \lambda^{|\sigma|} (\dd (u_{1},\ldots,u_{|\sigma|})) + \\
\nonumber \nonumber  \sum_{k=1}^{d} \binom{d}{k} \int_{(Y^{|\sigma|})_{[k]}} \left( \otimes_{j=0}^{d-1}  \left(\left( \bar{\mathbb{H}}^{d-j} \right)^{\otimes m_{j}}  \right) \right)_{\sigma} (u_{1},\ldots,u_{|\sigma|}) \times \\
\nonumber \left\vert \frac{ \frac{(d-k)!}{ (c-1-k)!}}{\frac{d!}{(c-1)!}}-\rho_{|\sigma|}(u_{1},\ldots,u_{|\sigma|}) \right\vert \lambda^{|\sigma|} (\dd (u_{1},\ldots,u_{|\sigma|})) \leq \\
\nonumber \sum_{k=1}^{c-1} \binom{c-1}{k} \int_{(Y^{|\sigma|})_{[k]}} \left( \otimes_{j=0}^{d-1}  \left(\left( \bar{\mathbb{H}}^{d-j} \right)^{\otimes m_{j}}  \right) \right)_{\sigma} (u_{1},\ldots,u_{|\sigma|}) \lambda^{|\sigma|} (\dd (u_{1},\ldots,u_{|\sigma|})) +\\
\nonumber  \sum_{k=1}^{d} \binom{d}{k} \int_{(Y^{|\sigma|})_{[k]}} \left( \otimes_{j=0}^{d-1}  \left(\left( \bar{\mathbb{H}}^{d-j} \right)^{\otimes m_{j}}  \right) \right)_{\sigma} (u_{1},\ldots,u_{|\sigma|}) \lambda^{|\sigma|} (\dd (u_{1},\ldots,u_{|\sigma|}))Re^{-Sa}
\end{gather}
where $(Y^{|\sigma|})_{[k]}$ is subspace of $Y^{|\sigma|}$, where facets $u_{1},\ldots,u_{|\sigma|}$ use orientations $e_{1},\ldots,e_{k}$ (each orientation is used at least by one of the fasets),$\binom{d}{k}$ is the number of possible selections of orientations used. We have an upper bound for the expression in the absolute value and we can get a lower bound in the same way.
 \hfill $\Box $
\begin{remark}
This means that we do not need the correlation function in the calculations and we can omit it from the formula, if we also omit $d-c+1$ possible orientations.
\end{remark}

\noindent{\bf Proof of Theorem}\\
It holds \cite{RefB}
\begin{align}
\label{Expect} & \E G_{j}(\mu^{(c)}_{a}) = \frac{a^{j}}{j!} \int_{Y^{j}} \mathbb{H}^{d-j}(\cap_{i=1}^{j} u_{i}) \rho_{j} (u_{1},\ldots,u_{j},\mu^{(c)}_{a}) \lambda^{j} (\dd (u_{1},\ldots,u_{j})),\\
\label{ProdExpect} & \E \prod_{j=1}^{c-1} G^{m_{j}}_{j}(\mu^{(c)}_{a}) = \sum_{\sigma \in \Pi_{1,\ldots,c-1}^{(m_{1},\ldots,m_{c-1})}} \prod_{j=1}^{c-1} \frac{1}{j!^{m_{j}}} a^{|\sigma|} \\
\nonumber & \int_{Y^{|\sigma|}} \left( \otimes_{j=1}^{c-1} \left(\left( \bar{\mathbb{H}}^{d-j} \right)^{\otimes m_{j}}  \right)  \right)_{\sigma}   (u_{1},\ldots,u_{|\sigma|}) \rho_{|\sigma|} (u_{1},\ldots,u_{|\sigma|},\mu^{(c)}_{a}) \lambda^{|\sigma|} (\dd (u_{1},\ldots,u_{|\sigma|}))  .
\end{align}
We can also get a relation for joint moments of centered random variables
\begin{align}
\label{CenProdExpect}  \E \prod_{j=1}^{c-1} \tilde{G}^{m_{j}}_{j}(\mu^{(c)}_{a}) = \frac{1}{a^{M}} \E \prod_{j=1}^{c-1} \left( G_{j}(\mu^{(c)}_{a})- \E G_{j}(\mu^{(c)}_{a}) \right)^{m_{j}} = \\
\nonumber \frac{1}{a^{M}} \sum_{i_{1}=0}^{m_{1}} \ldots \sum_{i_{c-1}=0}^{m_{c-1}} \binom{m_{1}}{i_{1}}\ldots \binom{m_{c-1}}{i_{c-1}} (-1)^{\sum_{j=1}^{c-1}i_{j}} \times \\
\nonumber \times \E \left( \prod_{j=1}^{c-1} G^{m_{j}-i_{j}}_{j}(\mu^{(c)}_{a}) \right) \prod_{j=1}^{c-1} \E \left( G_{j}(\mu^{(c)}_{a}) \right)^{i_{j}},\end{align}
where $\sum_{j=1}^{c-1}(j-\frac{1}{2})m_{j}=M$.\\
Firstly we calculate expectations of the U-statistics using Lemma \ref{SimpleLemma} and 
\begin{align}
\nonumber \frac{\E G_{j}(\mu^{(c)}_{a})}{a^{j}} & \rightarrow \frac{1}{j!}  \int_{Y^{j}_{c-1}} \mathbb{H}^{d-j} \left( \cap_{i=1}^{j} u_{i} \right) \lambda^{j}(\dd (u_{1},\ldots , u_{j}))= \\
\nonumber &= \frac{1}{d^{j}}  I_{j}\binom{c-1}{j},
\end{align}
where $\binom{c-1}{j}$ is number of possibilities how to select unique $j$ orientations from $c-1$ and $j!$ is number of possible allocation of them on the $j$ positions, $d^j$ is number of all possibilities how to select orientations from $d$ possibilities. $I_{j}$ is value of integrated Hausdorff measure over the space of centres, which does not depend on the currently selected orientations, they only need to be distinct, otherwise the measure would be 0, i.e. only non-paralel facets intersect and thus have non-zero volume of the intersection. We can already see that expectation of $U$-statistics of order higher or equal than $c$ are zero, therefore they converge to zero in $L^{2}$ and thus we only need to investigate the $U$-statistics of the order lower than $c$.\\
Secondly we calculate all joint moments. To do this we need to first use formula (\ref{CenProdExpect}) and Lemma \ref{SimpleLemma} - we use the limit values of correlation function, which we justify later. To describe the relation between original formula and the formula with correlation function replaced by its limit value we use $\simeq$
\begin{align}
\nonumber & \left( \prod_{j=1}^{c-1} j!^{m_{j}}\right) \E \left( \prod_{j=1}^{c-1} G^{m_{j}-i_{j}}_{j}(\mu^{(c)}_{a}) \right)\prod_{j=1}^{c-1} \E \left( G_{j}(\mu^{(c)}_{a}) \right)^{i_{j}} \simeq  \\
\nonumber & \simeq \prod_{j=1}^{c-1} \left( \int_{Y_{c-1}^{j}} \mathbb{H}^{d-j}(\cap_{i=1}^{j} u_{i}) \lambda^{j} (\dd (u_{1},\ldots,u_{j})) \right)^{i_{j}} \sum_{\sigma \in \Pi_{1,\ldots,c-1}^{(m_{1}-i_{1},\ldots,m_{c-1}-i_{c-1})}}  a ^{|\sigma|+\sum_{j=1}^{c-1}ji_{j}} \\
\nonumber & \int_{Y_{c-1}^{|\sigma|}} \left( \otimes_{j=1}^{c-1} \left( \left( \bar{\mathbb{H}}^{d-j} \right)^{\otimes (m_{j}-i_{j})}   \right)  \right)_{\sigma} (u_{1},\ldots,u_{|\sigma|}) \lambda^{|\sigma|} (\dd (u_{1},\ldots,u_{|\sigma|})) 
\end{align}

 We are interested only in terms with power higher than or equal to $M$, because the other terms will tend to zero with increasing $a$, i.e. partitions fullfilling condition $|\sigma| \geq M-\sum_{j=1}^{c-1}i_{j}j$. Also we do not have to examine odd moments, i.e. those with $\sum_{j=1}^{c-1}m_{j}$ odd, because there is not any summand with the power of $a$ matching $M$ in the divisor, thus they can be only zero or infinite, therefore if we prove that all even moments tend to some finite value, then all odd moments are equal to zero.\\
Select $\mathbf{s} = (s_{1},\ldots,s_{c-1} )$, so that $m_{i} \geq s_{i} \geq 0, i \in [c-1], \exists j \in [c-1], m_{j}>s_{j}$, choose any partition $\sigma_{\mathbf{s}} \in \Pi_{1,\ldots,c-1}^{\mathbf{s}}$ fullfilling conditions $|\sigma_{\mathbf{s}}| \geq M-\sum_{j=1}^{c-1}i_{j}j$ and $S(\sigma_{\mathbf{s}} )=0$, i.e. each block of $\pi$ is connected to any other block of $\pi$ by some block of $\mathbf{s}$. Then for $\mathbf{t}=(t_{1},\ldots,t_{c-1}), m_{i} \geq t_{i} \geq s_{i},i \in [c-1], \exists j \in [c-1],t_{j}>s_{j}$ there are partitions $\sigma_{\mathbf{t}} \in \Pi_{1,\ldots,c-1}^{\mathbf{t}}$, which have only additional singleton rows compared to $\sigma_{\mathbf{s}}$, $S(\sigma_{\mathbf{t}})=\sum_{i=1}^{c-1}t_{i}-s_{i}=|\sigma_{\mathbf{t}}|-|\sigma_{\mathbf{s}}|$ and it holds
\begin{gather}
\nonumber a^{|\sigma_{\mathbf{t}}|} \int_{Y^{|\sigma_{\mathbf{t}}|}_{c-1}} \left( \otimes_{j=1}^{c-1} \left( \left( \bar{\mathbb{H}}^{d-j} \right)^{\otimes t_{j}} \right) \right)_{\sigma_{\mathbf{t}}}(u_{1},\ldots,u_{|\sigma_{\mathbf{t}}|}) \lambda^{|\sigma_{\mathbf{t}}|} (\dd (u_{1},\ldots,u_{|\sigma_{\mathbf{t}}|})) = \\
\nonumber  a^{|\sigma_{\mathbf{s}}|} \int_{Y^{|\sigma_{\mathbf{s}}|}_{c-1}} \left( \otimes_{j=1}^{c-1} \left( \left( \bar{\mathbb{H}}^{d-j} \right)^{\otimes s_{j}} \right) \right)_{\sigma_{\mathbf{s}}}(u_{1},\ldots,u_{|\sigma_{\mathbf{s}}|}) \lambda^{|\sigma_{\mathbf{s}}|} (\dd (u_{1},\ldots,u_{|\sigma_{\mathbf{s}}|})) \times \\
\nonumber \times  \prod_{j=1}^{c-1}  \left( \left( \frac{a(c-1)}{d} \right)^{j} \int_{Y^{j}_{c-1}}\mathbb{H}^{d-j} \left( \cap_{i=1}^{j} u_{i} \right) \lambda^{j}(\dd (u_{1},\ldots , u_{j})) \right)^{t_{j}-s_{j}} ,
\end{gather}
because we can separate the singleton rows corresponding to the functions  $\bar{\mathbb{H}}^{k}$ in tensor product, which can be integrated separately, because they do not have any common variables with the other functions in the tensor product and the integral is equal to the expectation of $U$-statistic. We can see that all summands corresponding to any of the partitions $\sigma_{\mathbf{t}}$ in the evaluation of  (\ref{CenProdExpect}) contain common term
 \begin{gather}
\nonumber \Theta =  a^{|\sigma_{\mathbf{s}}|} \int_{Y^{|\sigma_{\mathbf{s}}|}_{c-1}} \left( \otimes_{j=1}^{c-1} \left( \left( \bar{\mathbb{H}}^{d-j} \right)^{\otimes s_{j}} \right) \right)_{\sigma_{\mathbf{s}}}(u_{1},\ldots,u_{|\sigma_{\mathbf{s}}|}) \lambda^{|\sigma_{\mathbf{s}}|} (\dd (u_{1},\ldots,u_{|\sigma_{\mathbf{s}}|})) \times \\
\nonumber \times  \prod_{j=1}^{c-1}  \left( a^{j} \int_{Y^{j}_{c-1}}\mathbb{H}^{d-j} \left( \cap_{i=1}^{j} u_{i} \right) \lambda^{j}(\dd (u_{1},\ldots , u_{j})) \right)^{m_{j}-s_{j}}
\end{gather}
and then we sum over all such partitions $\sigma_{\mathbf{t}}$
\begin{gather}
\nonumber \Theta \sum_{i_{1}=s_{1}}^{m_{1}} \ldots \sum_{i_{c-1}=s_{c-1}}^{m_{c-1}} \binom{m_{1}}{i_{1}}\ldots\binom{m_{c-1}}{i_{c-1}} \binom{i_{1}}{s_{1}}\ldots\binom{i_{c-1}}{s_{c-1}} (-1)^{\sum_{j=1}^{c-1}i_{j}}=\\
\nonumber \Theta (-1)^{\sum_{j=1}^{c-1}s_{j}} \binom{m_{1}}{s_{1}} \ldots \binom{m_{c-1}}{s_{c-1}} \times \\
\nonumber \times \sum_{i_{1}=0}^{m_{1}-s_{1}} \ldots \sum_{i_{c-1}=0}^{m_{c-1}-s_{c-1}} \binom{m_{1}-s_{1}}{i_{1}} \ldots \binom{m_{c-1}-s_{c-1}}{i_{c-1}} (-1)^{\sum_{j=1}^{c-1}i_{j}}=0,
\end{gather}
where we use Binomial theorem for summing with necessary condition $\sum_{j=1}^{c-1}s_{j}<\sum_{j=1}^{c-1}m_{j}$ and $\binom{m_{j}}{i_{j}}$ are original coefficients from formula (\ref{CenProdExpect}) and $\binom{i_{j}}{s_{j}}$ is the number of options how to select additional singleton rows. \\
Therefore all partitions with any singleton rows or containted within $\Pi^{\mathbf{s}}_{1,\ldots,c-1}, \mathbf{s} < \mathbf{m}=(m_{1},\ldots,m_{c-1})$ cancel each other out.
But we calculated with the limit values of correlation functions and the integrals are multiplied by $a$ in polynomial, thus we have to deal with speed of convergence. We have the limit in a form 
\begin{gather}
\nonumber a^{k} \left\vert \sum_{i=0}^{N} \int_{X} H(\mathbf{x})\varrho_{i}(\mathbf{x},a)\Lambda(\mathbf{x}) \right\vert \leq a^{k} \left\vert \sum_{i=0}^{N} \int_{X} H(\mathbf{x})|\varrho_{i}(\mathbf{x},a)-\varrho_{i}(\mathbf{x})|\Lambda(\mathbf{x}) \right\vert \leq  \\
\nonumber \leq a^{k} Re^{-Sa}  \left\vert \sum_{i=0}^{N} \int_{X} H(\mathbf{x})\Lambda(\mathbf{x}) \right\vert \rightarrow 0,
\end{gather}
where $H\geq 0$ , $\varrho_{i}\geq 0$,$| \varrho_{i}(\mathbf{x},a)-\varrho_{i}(\mathbf{x})|< Re^{-Sa}$ ,  $R,S>0$, \\$\sum_{i=0}^{N} \int_{X} H(\mathbf{x})\varrho_{i}(\mathbf{x})\Lambda(\mathbf{x})=0$, the sum is over all partitions, which are the same after leaving out all singleton rows and the measure represents product measure of $\lambda$ on space $X=Y^{p}$ for some $p$ and $K$ is created by multiplying Hausdorff measures and $\varrho_{i}$ by multiplying correlation functions.\\
Now we are left out only with partitions $\sigma$, which do not contain any pure singleton rows, are contained in $\Pi^{\mathbf{m}}_{1,\ldots,c-1}$. These partitions have each row connected exactly to one another row by one block of two elements in $\sigma$ ($|\sigma|=M$) and therefore, if we omit all the mentioned partitions, then \\  $ \left( \prod_{j=1}^{c-1} j!^{m_{j}}\right) \E \left( \prod_{j=1}^{c-1} G^{m_{j}-i_{j}}_{j}(\mu^{(c)}_{a}) \right)\prod_{j=1}^{c-1} \E \left( G_{j}(\mu^{(c)}_{a}) \right)^{i_{j}} \simeq$
\begin{gather}
\nonumber a^{M}\sum_{k^{(2)}_{1},\ldots,k^{(2)}_{m_{2}}=1}^{2}\ldots \sum_{k^{(c-1)}_{1},\ldots,k^{(c-1)}_{m_{c-1}}=1}^{c-1} \sum_{\sigma \in \tilde{\Pi}_{K},J\in \sigma: |J|=2} \prod_{J=\{b_{1},b_{2} \} \in \sigma} a^{\tau(b_{1}) + \tau(b_{2})-1} \\
\nonumber   \int_{Y_{c-1}^{\tau(b_{1}) + \tau(b_{2}) -1}} \mathbb{H}^{d-\tau(b_{1})}(\cap_{i=1}^{\tau(b_{1})}x_{i})\mathbb{H}^{d-\tau(b_{2})}(\cap_{i=1}^{\tau(b_{2})-1}x_{\tau(b_{1})+i} \cap x_{1})  \times \\
\nonumber \times \lambda^{\tau(b_{1}) + \tau(b_{2}) -1}(\dd (x_{1},\ldots,x_{\tau(b_{1}) + \tau(b_{2}) -1})),\tau(s)=\mathrm{max}_{j \in [c-1]} \left\{ \sum_{i=1}^{j-1}m_{i}<s  \right\}, K=\sum_{j=1}^{c-1}m_{j}
\end{gather}
where we sum first over all possible selections of common elements among the partitions and then over all possible pairings of partition rows, we also divide integral into several parts, where each part consists only of elements which are in the same block of a  partition. Function $\tau$ connects row of partition to its length.
It holds
\begin{gather}
\nonumber \int_{Y_{c-1}^{\tau(b_{1}) + \tau(b_{2}) -1}} \mathbb{H}^{d-\tau(b_{1})}(\cap_{i=1}^{\tau(b_{1})}x_{i}) \times \\
\nonumber \times \mathbb{H}^{d-\tau(b_{2})}(\cap_{i=1}^{\tau(b_{2})-1}x_{\tau(b_{1})+i} \cap x_{1})  
 \lambda^{\tau(b_{1}) + \tau(b_{2}) -1}(\dd (x_{1},\ldots,x_{\tau(b_{1}) + \tau(b_{2}) -1})) \\
\label{LastExp}=     \frac{(c-1)(\tau(b_{1})-1)!(\tau(b_{2})-1)! I_{\tau(b_{1})\tau(b_{2})} \binom{c-2}{\tau(b_{1})-1}\binom{c-2}{\tau(b_{2})-1}}{d^{\tau(b_{1})+\tau(b_{2})-1}} ,
\end{gather}
where $c-1$ is number of possibilities how to select the one common facet orientation, $\binom{c-2}{\tau(b_{1})-1}$, $\binom{c-2}{\tau(b_{2})-1}$ is number of possibilities how to select distinct remaining orientations of the rest of the facets in the first and the second function in integrand and $(\tau(b_{1})-1)!$,$(\tau(b_{2})-1)!$  are numbers of their possible positions among facets in the Hausdorff measure, $d^{\tau(b_{1})+\tau(b_{2})-1}$ is number of all possible orientations from all $d$ possibilities (even non-distinct ones) and $I_{\tau(b_{1})\tau(b_{2})}$ is integral over facets with fixed orientations over the space of the facet centres. Then  using (\ref{LastExp})\\$  \E \left( \prod_{j=1}^{c-1} G^{m_{j}-i_{j}}_{j}(\mu^{(c)}_{a}) \right)\prod_{j=1}^{c-1} \E \left( G_{j}(\mu^{(c)}_{a}) \right)^{i_{j}} \simeq$
\begin{gather}
\nonumber  \left(\frac{a}{d}\right)^{M} \sum_{\sigma \in \tilde{\Pi}_{K},J\in \sigma: |J|=2} \sum_{k^{(2)}_{1},\ldots,k^{(2)}_{m_{2}}=1}^{2}\ldots \sum_{k^{(c-1)}_{1},\ldots,k^{(c-1)}_{m_{c-1}}=1}^{c-1} \prod_{J=\{b_{1},b_{2} \} \in \sigma} \\
\nonumber \left( \prod_{j=1}^{c-1}\frac{1}{j!^{m_{j}}} \right) (\tau(b_{1})-1)!(\tau(b_{2})-1)!(c-1)I_{\tau(b_{1})\tau(b_{2})} \binom{c-2}{\tau(b_{1})-1}\binom{c-2}{\tau(b_{2})-1}=\\
\nonumber \left( \frac{a}{d} \right)^{M}\sum_{\sigma \in \tilde{\Pi}_{K},J\in \sigma: |J|=2} \prod_{J=\{b_{1},b_{2} \} \in \sigma} (c-1)I_{\tau(b_{1})\tau(b_{2})} \binom{c-2}{\tau(b_{1})-1}\binom{c-2}{\tau(b_{2})-1}.
\end{gather}
If we express the covariance of any two variables in the same form by selecting $\tilde{G}_{i}(\mu^{(c)}_{a})$ and $\tilde{G}_{j}(\mu^{(c)}_{a}) $ and then by calculating
\begin{align}
\nonumber \frac{a^{i+j-1}(c-1)I_{ij}}{d^{i+j-1}} \binom{c-2}{i-1}\binom{c-2}{j-1}
\end{align}
we can see that the distribution of statistics has the property of normal distribution, i.e. joint moments of centered variables are equal to sum over all pairs of unordered random variables (random variables with higher power are used as several distinct multiplied random variables) and this implies the central limit theorem, because normal distribution is defined by its moments \cite{RefBil}.\\
There is only one remaining statement to prove
\begin{align}
\nonumber &\frac{G_{j}(\mu^{(c)}_{a})}{a^{j}} \xrightarrow{L^{2}} \frac{I_{j}}{d^{j}} \binom{c-1}{j}, ~ & c \in \{ 2,\ldots,d \}, j < c ,
\end{align}
the first moment of the random variable on the left side is equal to right side and the variance tends to zero as can be seen from the central limit theorem.
 \hfill $\Box $

\begin{remark}
Consider process $\mu_{a}$ with density in more general form $p(\mathbf{x})=c_{a} \exp \left(  \sum_{i=1}^{d} \nu_{i} G_{i}(\mathbf{x}) \right)$. Assume there is $c \geq 2, \nu_{c}>0$ and select such minimal $c$. 
\begin{align}
\nonumber &\E \exp \left(  \sum_{j=1}^{d} \nu_{j} G_{j}(\eta_{a}) \right) = \\
\nonumber &= \sum_{n=0}^{\infty}\frac{a^{n} e^{-aT}}{n!} \int_{Y^{n}} \exp \left(  \sum_{j=1}^{d} \nu_{j} G_{j}(\{ u_{1},\ldots,u_{n} \}) \right) \lambda^{n} (\dd (u_{1},\ldots,u_{n})) \\
\nonumber  &\geq e^{-aT} \sum_{n_{1}=0}^{\infty} \ldots  \sum_{n_{d}=0}^{\infty} \frac{a^{n_{1}+\ldots+n_{d}}}{n_{1}!\ldots n_{d}!} \exp \left( \sum_{j=1}^{d} \nu^{\prime}_{j} \sum_{\{i_{1},\ldots,i_{j} \}\subset [d] } \prod_{l=1}^{j}n_{i_{l}}
 \right) \\
\nonumber  &\geq e^{-aT} \sum_{n_{1}=0}^{\infty} \ldots  \sum_{n_{c}=0}^{\infty} \frac{a^{n_{1}+\ldots+n_{c}}}{n_{1}!\ldots n_{c}!} \exp \left( \sum_{j=1}^{c} \nu^{\prime}_{j} \sum_{\{i_{1},\ldots,i_{j} \}\subset [c] } \prod_{l=1}^{j}n_{i_{l}}
 \right) \\
 \label{divSuma} &\geq e^{-aT} \sum_{n=0}^{\infty} \frac{a^{nc}}{(n!)^{c}} \exp \left( 
 \sum_{j=1}^{c} \nu^{\prime}_{j} \binom{c}{j} n^{j} \right),
\end{align}
where $\nu^{\prime}_{j}\nu_{j}>0$ and we firstly set the last $d-c$ summing variables to zero and then summed only over the summands, where all of the summing variables have the same value. It can be proven (e.g. by using ratio test), that the sum in ($\ref{divSuma}$) is divergent and therefore $p(\eta_{a}) \not\in L^{1}(P_{\eta_{a}})$ in this case. On the other hand non-positivity of parameters $\nu$ implies $p(\eta_{a}) \in L^{1}(P_{\eta_{a}}) \cap L^{2}(P_{\eta_{a}})$ as shown in \cite{RefB}, which finally leads to  $\nu_{l} \leq 0, l \geq 2 \Longleftrightarrow p(\eta_{a}) \in L^{1}(P_{\eta_{a}}) \cap L^{2}(P_{\eta_{a}})$.
\end{remark}

\begin{remark}
Consider process $\mu_{a}$ with density $p(\mathbf{x})=c_{a} \exp \left(  \sum_{i=1}^{d} \nu_{i} G_{i}(\mathbf{x}) \right)$, $\nu_{l} \leq 0, l \geq 2$. Assume there is $c \geq 2, \nu_{c}<0$ and select minimal such $c$. Then using similar techniques as in proof of Lemma \ref{KorelFce} and Lemma \ref{LemmaLim} we can show that 
\begin{align}
\nonumber \lim_{a\rightarrow \infty} \rho_{p}(x_{1},\ldots,x_{p},\mu_{a}) = \lim_{a\rightarrow \infty} \rho_{p}(x_{1},\ldots,x_{p},\mu^{(c)}_{a}),
\end{align}
which leads to the same asymptotic distribution of statistics $(\tilde{G}_{1}(\mu_{a}),\ldots,\tilde{G}_{d}(\mu_{a}))$ as $(\tilde{G}_{1}(\mu^{(c)}_{a}),\ldots,\tilde{G}_{d}(\mu^{(c)}_{a}))$.
\end{remark}

\paragraph{Acknowledgements}
This research was supported by grant SVV 260225 of Charles University in Prague.

\end{document}